\documentclass{article}

\usepackage{graphicx,xcolor}
\usepackage{amssymb,amsthm,amsmath}
\DeclareMathOperator{\D}{d}
\newtheorem{thm}{Theorem}[section]

\newtheorem{prop}[thm]{Proposition}

\theoremstyle{definition}

\theoremstyle{remark}

\usepackage[font=small,labelfont=bf,tableposition=top]{caption}
\usepackage[font=footnotesize]{subfig}
\usepackage[demo,abs]{overpic}

\newcommand{\la}{\lambda}
\newcommand{\re}{{\rm e}} 
\newcommand{\ri}{{\rm i}} 

\title{Linear dispersive shocks}

\author{David Smith$^1$, Thomas Trogdon$^2$, Vishal Vasan$^3$ \\ \footnotesize $1$ \texttt{dave.smith@yale-nus.edu.sg}, Yale-NUS College, Corresponding author, \\ \footnotesize $2$ \texttt{trogdon@uw.edu}, trogdon@uw.edu \\ \footnotesize $3$ \texttt{vishal.vasan@icts.res.in}, International Centre for Theoretical Science \\ \footnotesize and Tata Institute of Fundamental Research}

\begin{document}

\maketitle
\begin{abstract}
We present a linear dispersive partial differential equation which manifests a number of qualitative features of dispersive shocks, typically thought to occur only in nonlinear models. The model captures much of the short time phenomenon but deviates from the full nonlinear model in its long time behavior. Though we limit our present discussion to dispersive shocks occurring in the Korteweg-de Vries equation, our work extends readily to other higher order dispersive models too.
\end{abstract}

\section{Dispersive shocks in the Korteweg-de Vries equation}
\label{s:kdv}

We consider a model of the Korteweg-de Vries (KdV) equation
\begin{align}\label{e:kdv}
    u_t + 6 F(u;x,t) u_x + u_{xxx} = 0, \quad (x,t) \in \mathbb R \times (0,\infty).
\end{align}
with initial condition $u(x,0) = u_0(x)$ with step-like initial datum that satisfies
\begin{align*}
    \lim_{x \to -\infty} u_0(x) = a > 0, \quad \lim_{x \to +\infty} u_0(x) = 0.
\end{align*}
If one chooses $F(u;x,t) = u(x,t)$, then the problem can be effectively solved with the inverse scattering transform (IST) \cite{Ablowitz2013,Egorova2013,Bilman2018}. But we emphasize that the IST does not produce an explicit solution expression for the initial value problem.

Solutions to \eqref{e:kdv} with $F(u;x,t)=u(x,t)$ are characterized by a dispersive shock wave that moves with velocity $+ 4a$ and a finite number of solitons that move with velocity $> 4 a$ \cite{Ablowitz2013,Bilman2018}.  In this paper we introduce a simple alternative choice for $F(u;x,t)$ that produces, for short times, the same qualitative dispersive shock wave.  Importantly, our choice turns \eqref{e:kdv} into a linear interface problem \cite{Sheils2015b} and an explicit solution expression is available.

\section{A linear moving-interface problem}

We begin by replacing \eqref{e:kdv} by a suitable linear model. This amounts to making a specific choice for $F(u;x,t)$. Specifically, we choose
\begin{align}\label{eq:average}
    F(u;x,t) = \begin{cases}\lim_{R \to \infty} \frac{1}{R} \int_{ct-R}^{ct} u(s,t) \D s & x < ct\\
    \lim_{R \to \infty} \frac{1}{R} \int_{ct}^{ct+R} u(s,t) \D s & x \geq ct, \end{cases} \quad c > 0.
\end{align}
Here $c$ is the front speed and is arbitrary. For example, $c = 4a$ is an obvious choice to model the case $F(u;x,t) = u(x,t)$.  We leave it as a free parameter.  Then \eqref{e:kdv} becomes
\begin{align}\label{eq:model}
\begin{split}
    u_t + u_{xxx} &= - u_x\begin{cases} 6a & x < ct,\\
    0 & x > ct. \end{cases}, \\
    \lim_{x \to - \infty} u(x,t) &= a, \quad
    \lim_{x \to + \infty} u(x,t) = 0,\\
    u(x,0) &= u_0(x).
    \end{split}
\end{align}
The discontinuity in the coefficient of the $u_x$ term is interpreted by enforcing that the solution, and its first and second derivatives are continuous in $x$ at $x = ct$.

We immediately move into a traveling frame by defining $q(x,t) = u(x + ct, t)$.  Then because $q_t = c u_x + u_t$, $q_x = u_x$ we find
\begin{subequations} \label{e:q.problem}
\begin{align}
    q_t + q_{xxx} &= q_x\begin{cases} c-6a & x < 0,\\
    c & x > 0, \end{cases} \\
    \lim_{x \to - \infty} q(x,t) &= a, \quad
    \lim_{x \to + \infty} q(x,t) = 0,\\
    q(x,0) &= u_0(x).
\end{align}
In order to keep things fully explicit, we choose
\begin{align}\label{eq:stepdata}
    u_0(x) = \begin{cases}
    a & x < 0,\\
    0 & x \geq 0. \end{cases}
\end{align}
and have the interface (continuity) conditions
\begin{align}
\lim_{x \to 0^-}\partial_x^j q(x,t) = \lim_{x \to 0^+} \partial_x^j q(x,t), \quad j =0,1,2, t \geq 0.
\end{align}
\end{subequations}
Throughout this paper, $q(x,t) = q(x,t;a,c)$ will refer to this solution.
\section{Explicit solution via UTM}

The unified transform method (UTM) was adapted in \cite{Sheils2015b} to produce an integral expression for the solution of the time-dependent linear Schr\"{o}dinger equation with piecewise constant potential.
The below proposition represents a generalization of that method to problem~\eqref{e:q.problem}, which is of a higher spatial order but otherwise similar.
The solution formula presented below owes its relative simplicity to the fact that the transformed problem for $q(x,t)-a\chi_{(-\infty,0)}(x)$ has homogeneous initial condition, where $\chi_A$ is the characteristic function of the set $A$.

\begin{prop}\label{prop}
    Problem~\eqref{e:q.problem} has solution representation, valid for arbitrary $\tau\geq t$, for $x<0$,
    \begin{multline*}
        2\pi q(x,t) = 2\pi a {}+{} \\
        \int_{\partial D} \re^{\ri\nu_1(\la)x+\ri\la^3t} \nu'_1(\la) \bigg[ g_2(-\ri\la^3;\tau) + \ri \nu_1(\la) g_1(-\ri\la^3;\tau) \\
        {}-{} \left( [\nu_1(\la)]^2+c-6a \right) \left\{ g_0(-\ri\la^3;\tau) - ah(-\ri\la^3;\tau) \right\} \bigg] \D\la \\
        {}+{} \int_{\partial D} \re^{\ri\nu_2(\la)x+\ri\la^3t} \nu'_2(\la) \bigg[ g_2(-\ri\la^3;\tau) + \ri \nu_2(\la) g_1(-\ri\la^3;\tau) \\
        {}-{} \left( [\nu_2(\la)]^2+c-6a \right) \left\{ g_0(-\ri\la^3;\tau) - ah(-\ri\la^3;\tau) \right\} \bigg] \D\la,
    \end{multline*}
    and, for $x>0$,
    \begin{multline*}
        2\pi q(x,t) = \int_{\partial D} \re^{\ri\nu_0(\la)x+\ri\la^3t} \nu'_0(\la) \bigg[ g_2(-\ri\la^3;\tau) \\
        {}+{} \ri \nu_0(\la) g_1(-\ri\la^3;\tau) - \left( [\nu_0(\la)]^2+c \right) g_0(-\ri\la^3;\tau) \bigg] \D\la,
    \end{multline*}
    where, outside a bounded region containing all branch cuts, the maps $\nu_j$ are analytic and satisfy
    \begin{align*}
        \nu_0^3+c\nu_0&=\la^3, \\
        \nu_1(\la)&=\nu_2(\alpha^2\la), \\
        \nu_2^3+(c-6a)\nu_2&=\la^3,
    \end{align*}
    $\nu_j(\la)/\la \approx \alpha^j$ for $\lvert\la\rvert$ large.  Additionally,
    \begin{align*}
        h(\rho;\tau) &= \int_0^\tau \re^{\rho s} \D s = \frac{1}{\rho} \left(\re^{\rho\tau}-1\right), \\
        \alpha &= \exp(2\pi\ri/3),
    \end{align*}
    the functions $g_j$ satisfy the linear system
    \begin{equation*}
        \mathcal{A}(\la) \begin{pmatrix} g_2(-\ri\la^3;\tau) \\ \ri g_1(-\ri\la^3;\tau) \\ - g_0(-\ri\la^3;\tau) \end{pmatrix} 
        {}={} \begin{pmatrix} 0 \\ 0 \\ - a([nu_2(\alpha\la)]^2+c-6a)h(-\ri\la^3;\tau) \end{pmatrix},
    \end{equation*}
    for
    \begin{equation*}
        \mathcal{A}(\la) =
        \begin{pmatrix}
            1 & \nu_0(\alpha\la) & [\nu_0(\alpha\la)]^2 + c \\
            1 & \nu_0(\alpha^2\la) & [\nu_0(\alpha^2\la)]^2 + c \\
            1 & \nu_2(\alpha\la) & [\nu_2(\alpha\la)]^2 + c - 6a
        \end{pmatrix},
    \end{equation*}
    $D$ is the sector $\pi/3<\arg\la<2\pi/3$ excluding any choice of compact simply connected set such that all integrands are analytic on the interior of $D$ and admit continuous bounded extension to $\partial D$ (the positively oriented boundary of $D$).
\end{prop}

The proof of this proposition may be obtained by following the method of~\cite{Sheils2015b}, suitably modified in light of the results of~\cite{Deconinck}, after a further change of variables to subtract the initial datum from $q$, so that the new problem has inhomogeneity appearing only as a jump datum in the interface condition.

This solution representation in Proposition~\ref{prop} has the advantage that contours can be deformed in the complex plane to make use of regions where the exponentials decay.  Once such contours are found, our methodology employs standard quadrature routines to compute the solution.  The methodology we employ in this paper is not asymptotically accurate (i.e. accurate for all values of $(x,t)$) but it is more than sufficient to present high-accuracy plots for moderate values of $(x,t)$.

\section{Solution analysis}

Before we analyze the solution expression provided by Proposition \ref{prop} we recall some relevant facts. Firstly, the KdV equation, linearized about the solution $u(x,t) = 0$ with initial datum \eqref{eq:stepdata} has a similarity solution, see \cite{GinoTomIVP}, for example.  The same is not true in our setting \eqref{eq:average}. This is an important distinction between our linearization and a naive one.  Secondly, the solution profile for a similarity solution of the form
\begin{align*}
    u(x,t) = U(x/t^{1/3}), \quad t > 0,
\end{align*}
will have the same maximum amplitude for all $t > 0$. Lastly, we state three main characteristics of the initial-value problem \eqref{e:kdv}--\eqref{eq:average} with \eqref{eq:stepdata}:
\begin{itemize}
    \item The maximum of the solution profile grows in time, at least for small times.  This is a property of the KdV equation, stemming from its nonlinearity, that is preserved in our linear model \eqref{eq:model}.
    \item The increase in amplitude of initial datum (an increase in $a$ in \eqref{eq:stepdata}) results in higher frequency oscillations.
    \item Lastly, \eqref{eq:stepdata} gives the universal limit for general small amplitude data.
\end{itemize}
A simple re-scaling permits us to combine the two independent parameters $a$ and $c$ into a single dynamically relevant parameter. Indeed given $q(x,t;c,a)$, $c > 0$, $a > 0$, consider
\begin{align*}
    Q(x,t) = \frac{1}{a} q(x c^{-1/2},t c^{-3/2}).
\end{align*}
Then $Q(x,0) = \begin{cases} 1 & x < 0, \\ 0 & x \geq 0.\end{cases}$ and for $x < 0$
\begin{align*}
    0 = q_t + q_{xxx} - (c - 6a) q_x = {ac^{3/2}}\left( Q_t + Q_{xxx} - \frac{c - 6a}{c} Q_x \right)
\end{align*}
while for $x > 0$
\begin{align*}
   0 =  q_t + q_{xxx} - c q_x = {ac^{3/2}}\left( Q_t + Q_{xxx} - Q_x \right).
\end{align*}
From this we conclude that $Q(x,t) = q(x,t;1,a/c)$ and thus we limit the analysis to the function $Q(x,t)$.

\subsection{Short time}

\begin{figure}[tbp]
\begin{center}
\subfloat[A plot of $U(x,0.1)$.\label{f:t-0.1}]{
    \begin{overpic}[width=0.473\textwidth]{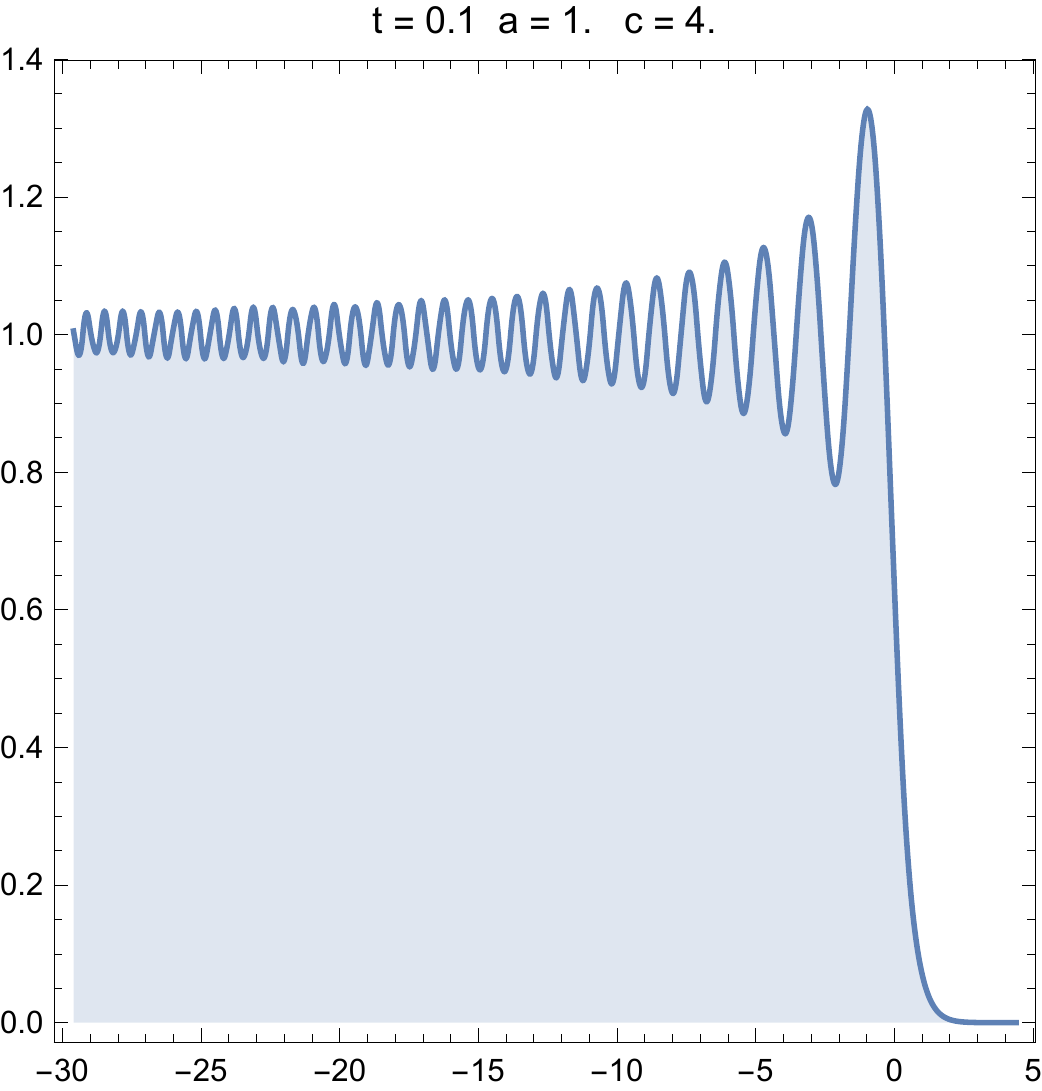}
        \put(45,166){\colorbox{white}{$t=0.1$ $a=1$ $c=4$}}
        \put(167,5){$x$}
    \end{overpic}
}  \\
\subfloat[A plot of $U(x,0.5)$.\label{f:t-0.5}]{
    \begin{overpic}[width=0.473\textwidth]{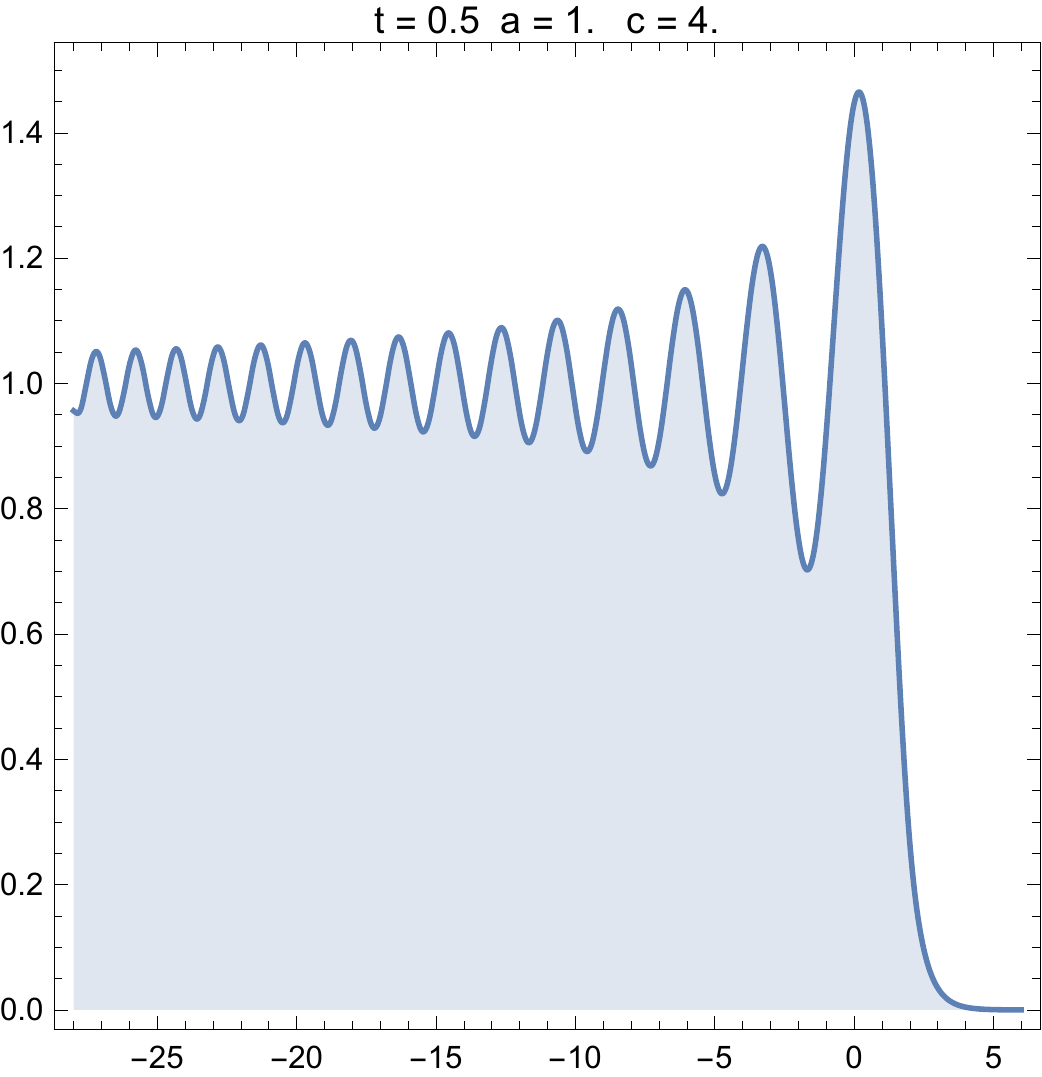}
        \put(45,166){\colorbox{white}{$t=0.5$ $a=1$ $c=4$}}
        \put(167,5){$x$}
    \end{overpic}
}
\end{center}
\caption{The function $U(x,t)$ for small times}
\end{figure}

\begin{figure}[tbp]
\begin{center}
\subfloat[A plot of $U(x,1.5)$.\label{f:t-1.5}]{
    \begin{overpic}[width=0.473\textwidth]{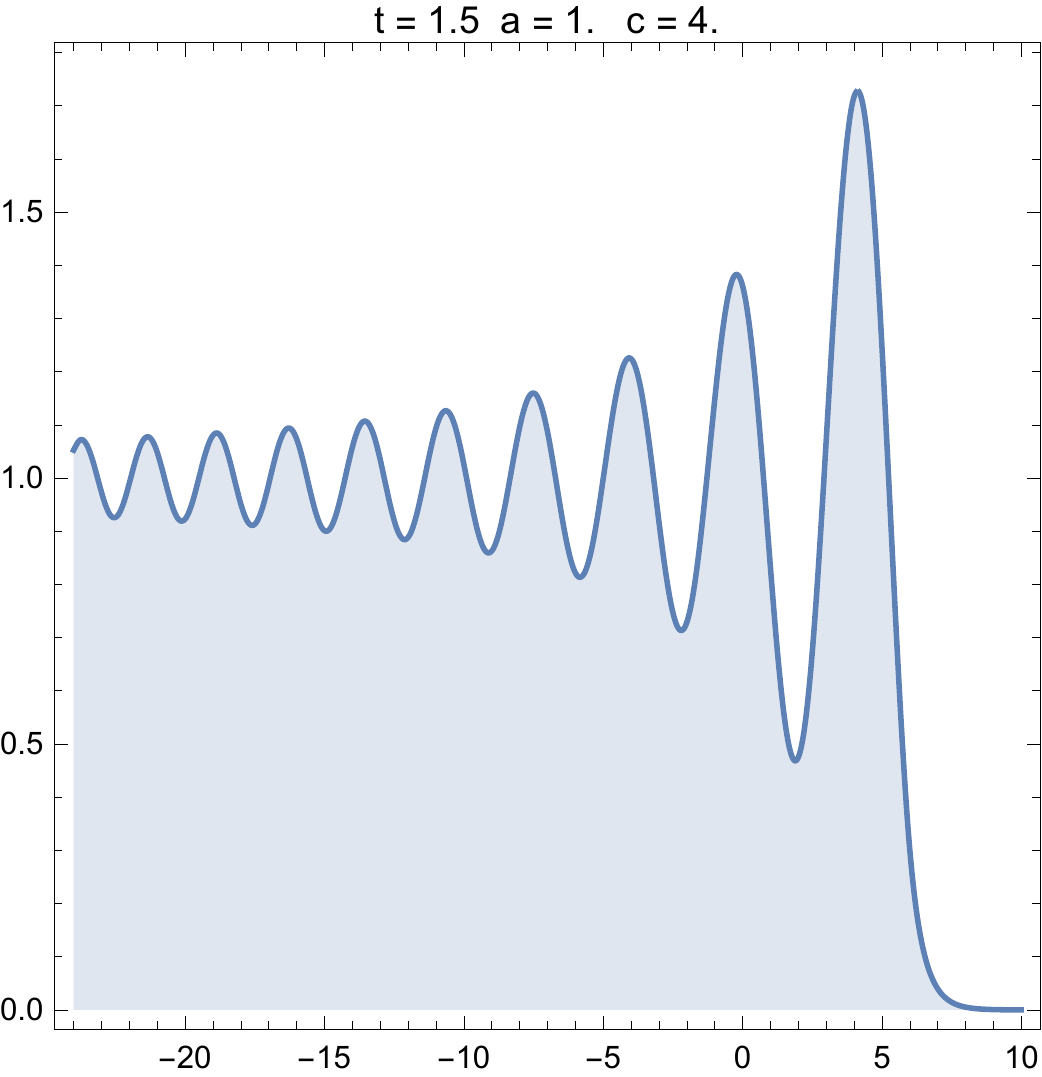}
        \put(45,166){\colorbox{white}{$t=1.5$ $a=1$ $c=4$}}
        \put(167,5){$x$}
    \end{overpic}
}  \\
\subfloat[A plot of $U(x,2.75)$ with a line to show the nearly linear arrangement of the first three peaks.\label{f:t-2.75}]{
    \begin{overpic}[width=0.473\textwidth]{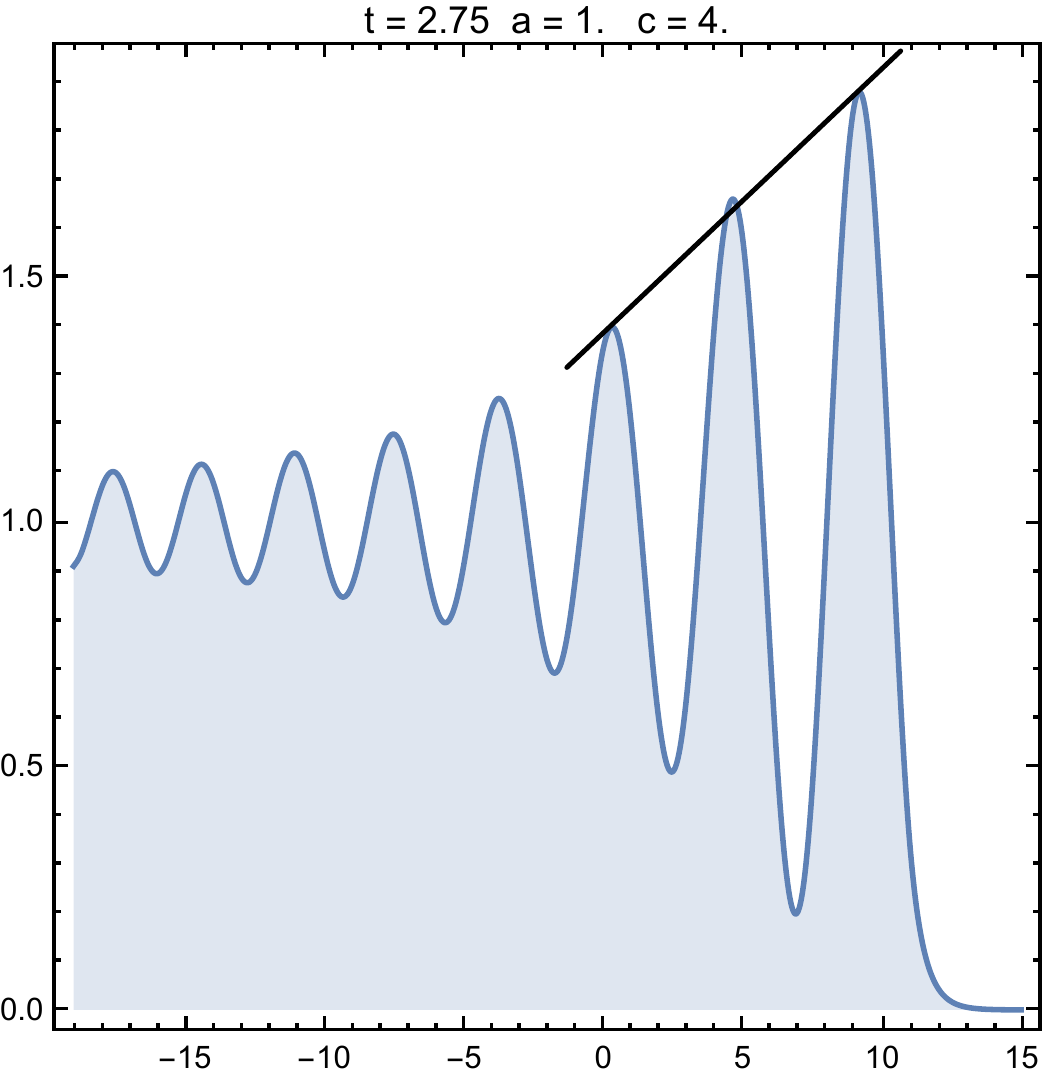}
        \put(42,166){\colorbox{white}{$t=2.75$ $a=1$ $c=4$}}
        \put(167,5){$x$}
    \end{overpic}
}
\end{center}
\caption{The function $U(x,t)$ for small, but larger times}
\end{figure}

We analyze the solution for short time and for large values of $a$ and $c$.  Define $\gamma = a/c$ and assume $\gamma > 1/6$ is fixed.
As $c$ increases, analyzing $U(x,t):= Q(x-\gamma t,t)$ at a fixed time amounts to analyzing $u(x,t)$ (the solution to \eqref{eq:model}) at a time $\propto c^{-3/2}$.  For this reason we consider this to be a small time analysis.  In Figure~\ref{f:t-0.1} we plot $U(x,0.1)$ to show the high oscillation that is produced by the initial incompatibility between initial condition and interface (continuity) conditions.  The amplitude of the solution increases as time increases in Figures~\ref{f:t-0.5}, \ref{f:t-1.5} and \ref{f:t-2.75}.

\subsection{Intermediate time}

\begin{figure}[tbp]
\begin{center}
\subfloat[A plot of $U(x,4)$.\label{f:t-4}]{
    \begin{overpic}[width=0.473\textwidth]{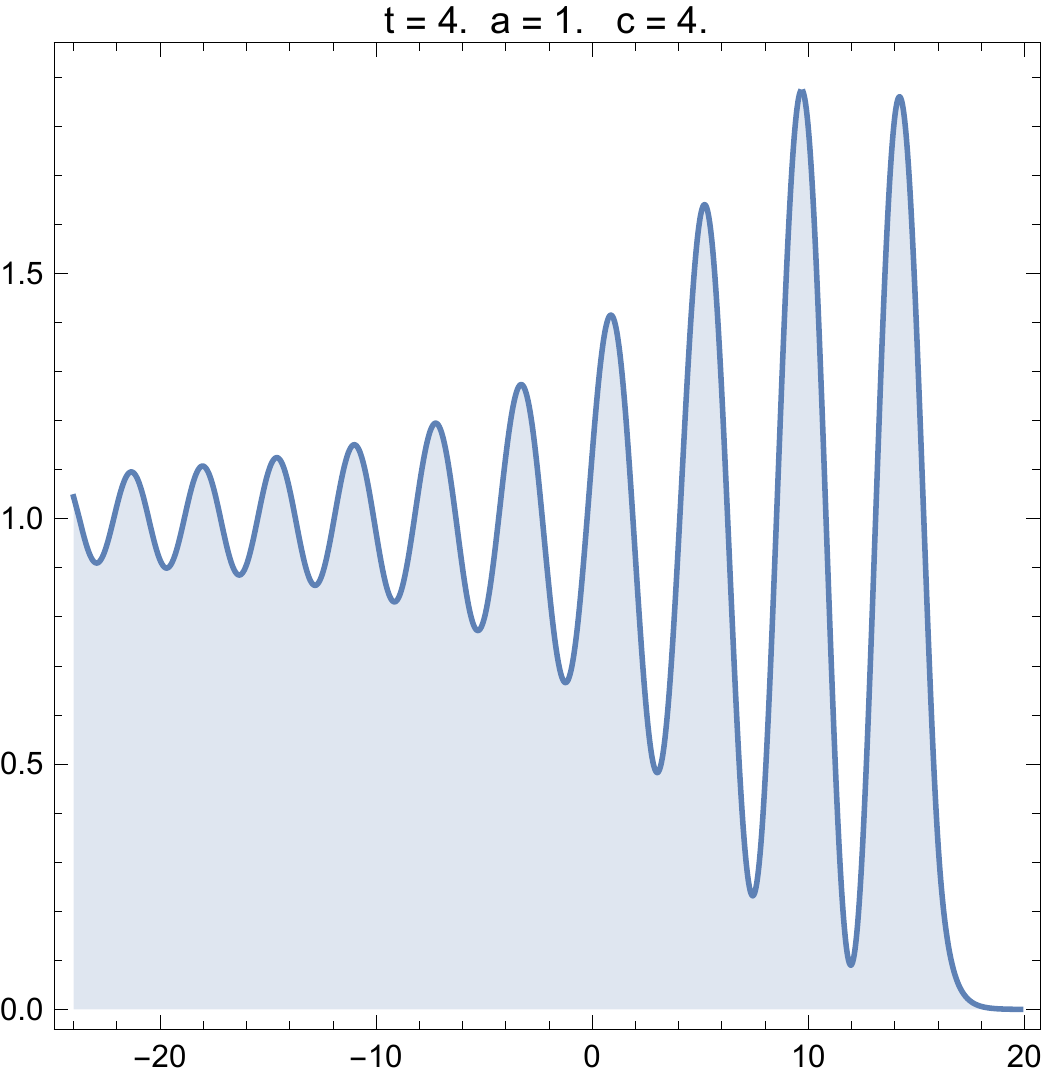}
        \put(47,166){\colorbox{white}{$t=4$ $a=1$ $c=4$}}
        \put(167,5){$x$}
    \end{overpic}
}  \\
\subfloat[A plot of $U(x,8)$.\label{f:t-8}]{
    \begin{overpic}[width=0.473\textwidth]{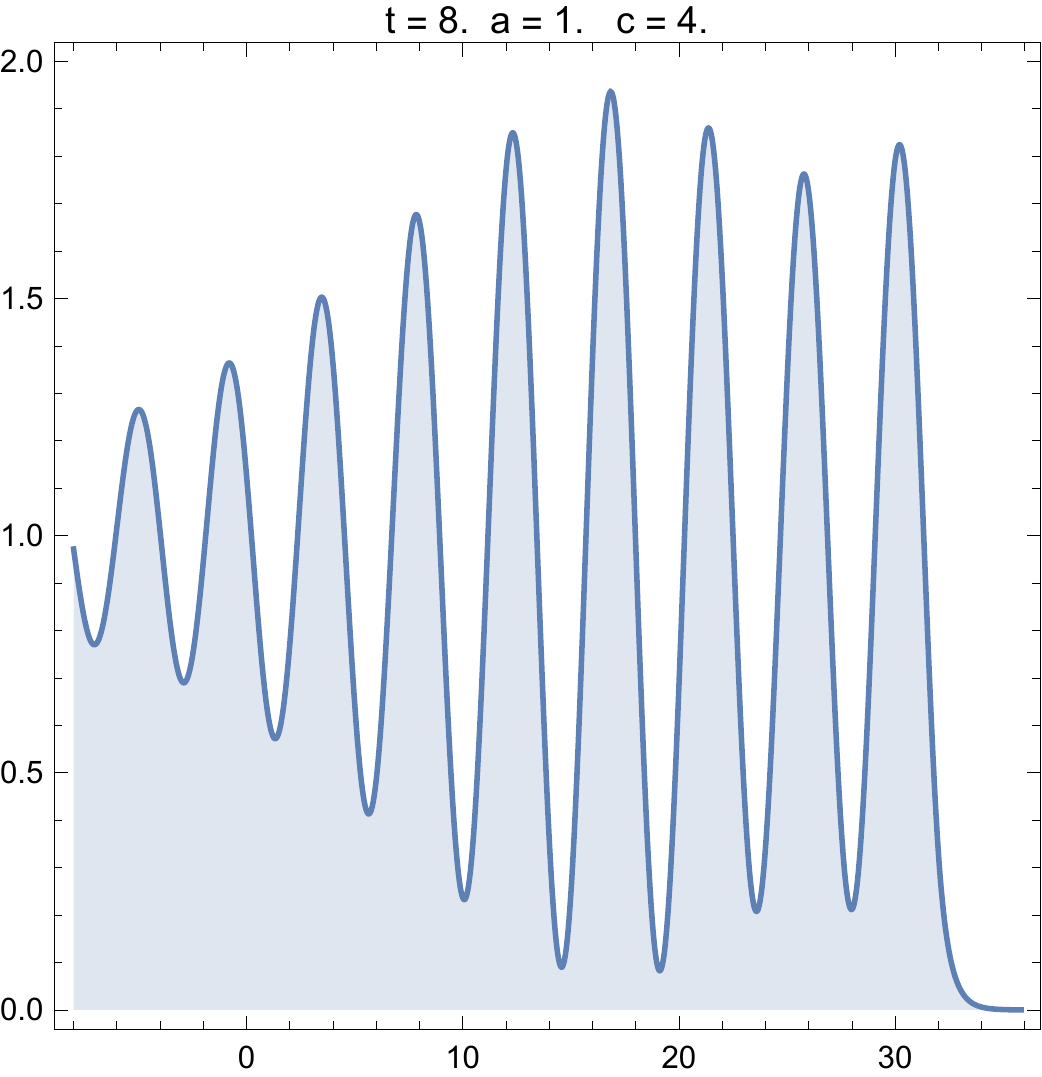}
        \put(47,166){\colorbox{white}{$t=8$ $a=1$ $c=4$}}
        \put(167,5){$x$}
    \end{overpic}
}
\end{center}
\caption{The function $U(x,t)$ at intermediate times showing a pulse, or wave packet, that is generated at the interface $x = ct$.}
\end{figure}

The linear arrangement of the peaks, reminiscent of a similar arrangement in dispersive shocks produced in KdV,  shown in Figure~\ref{f:t-2.75} is transient.  In Figures~\ref{f:t-4} and \ref{f:t-8} we show that the traveling interface generates a pulse that travels in the negative $x$-direction.  This pulse leaves lower amplitude oscillations in its wake.

\subsection{Long time}

\begin{figure}[tbp]
\begin{center}
    \begin{overpic}[width=0.473\textwidth]{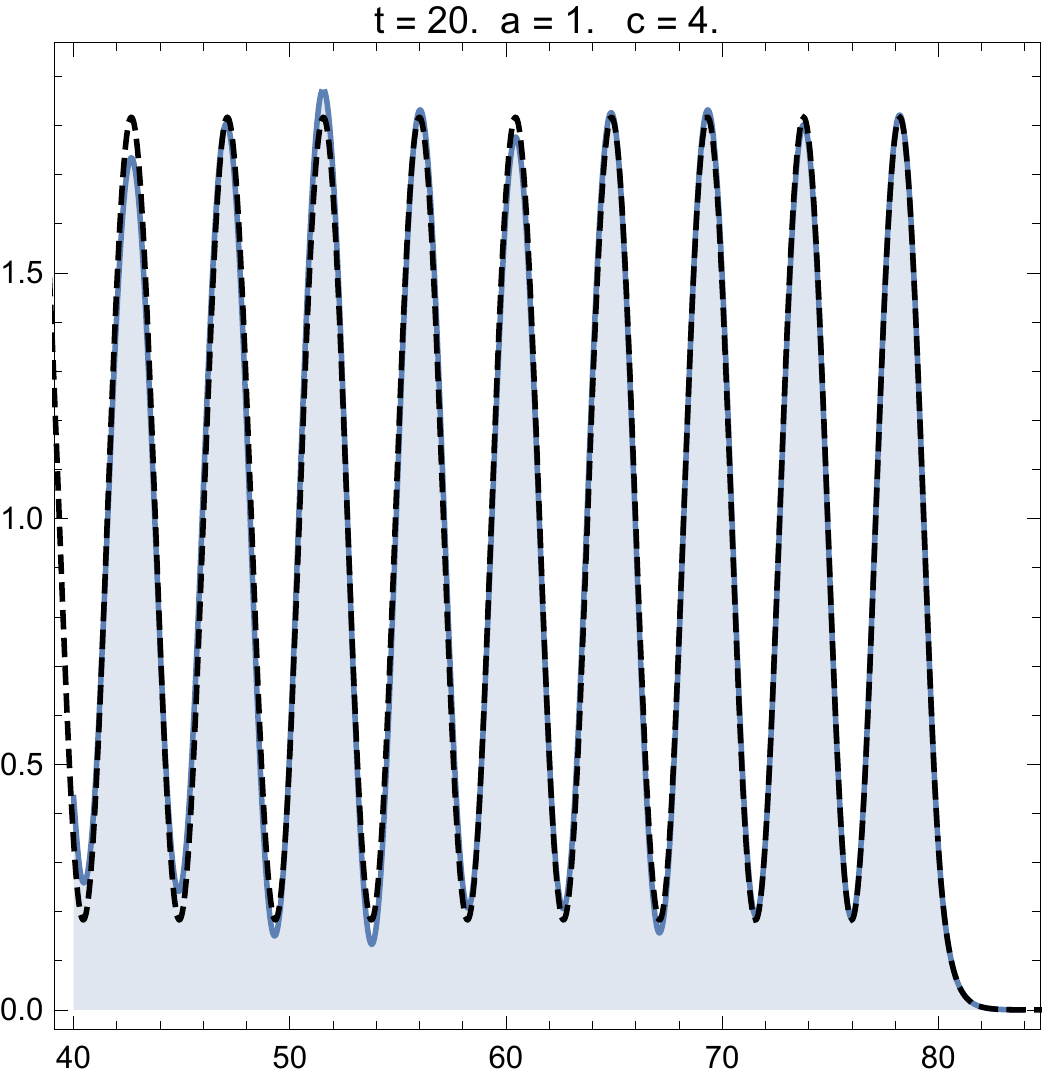}
        \put(45,166){\colorbox{white}{$t=20$ $a=1$ $c=4$}}
        \put(167,5){$x$}
    \end{overpic}
\end{center}
\caption{The function $U(x,t)$ at $t = 20$ to show how it approaches the long-time limit given by the black dashed curve. \label{f:t-20}}
\end{figure}

The long-time behavior of the solution is easily determined by looking for a stationary solution: Solve
\begin{align*}
     q_{xxx} &= c q_x, \quad x > 0,\\
     q_{xxx} &= (c - 6a) q_x, \quad x < 0.
\end{align*}
Assuming that $c - 6a < 0$ we find that the solution is given by
\begin{align*}
    q(x,t) \sim \begin{cases} c_1 + c_2 \re^{ \sqrt{c} x} + c_3 \re^{ -\sqrt{c} x} \hspace{5em}\mbox{if } x > 0,\\
    b_1 + b_2 \sin \left( x\sqrt{6a-c}\right) + b_3 \cos\left( x\sqrt{6a-c} \right) \\ \phantom{c_1 + c_2 \re^{ \sqrt{c} x} + c_3 \re^{ -\sqrt{c} x}}\hspace{5em}\mbox{if } x < 0.
    \end{cases}
\end{align*}
The boundary and interface conditions imply $c_1 = c_2 = 0$, $b_1 = a$, $c_3 = a - c/6$, $b_2 = -\sqrt{\frac{c}{6a-c}}(a - c/6)$ and $b_3 = -c/6$.  We compare this stationary solution with the solution $U(x,20)$ in Figure~\ref{f:t-20}.

\subsection{Small amplitude}

We point out that, after rescaling, $U(x,t)$ is the small-amplitude limit of $u(x,t)$ with general initial data.  Indeed let $u(x,0) = u_0(x) = a ( v(x) + h(x) )$ where
\begin{align*}
    h(x) = \begin{cases} 1 & x < 0,\\ 0 & x \geq 0 \end{cases}
\end{align*}
and $v(x)$ is an integrable function with, for simplicity, compact support.  Suppose again that $\gamma > 1/6$ is fixed.  Then, define
\begin{align*}
    \tilde U(x,t) = \frac{1}{a} u\left(x c^{-1/2}+c^{-1/2}t - \frac{a}{c^{-3/2}}t,c^{-3/2}t\right)
\end{align*}
It follows that $\tilde U$ is a solution of the same equation as $U$ but now with initial datum
\begin{align*}
    \tilde U(x,0) = v(x c^{-1/2}) + h(x).
\end{align*}
As $c \downarrow 0$, the initial datum converges to $h(x)$;  $L^1$ convergence of initial datum will imply uniform convergence on compact sets for $t > 0$.

\begin{figure}[tbp]
\begin{center}
    \begin{overpic}[width=0.56\textwidth]{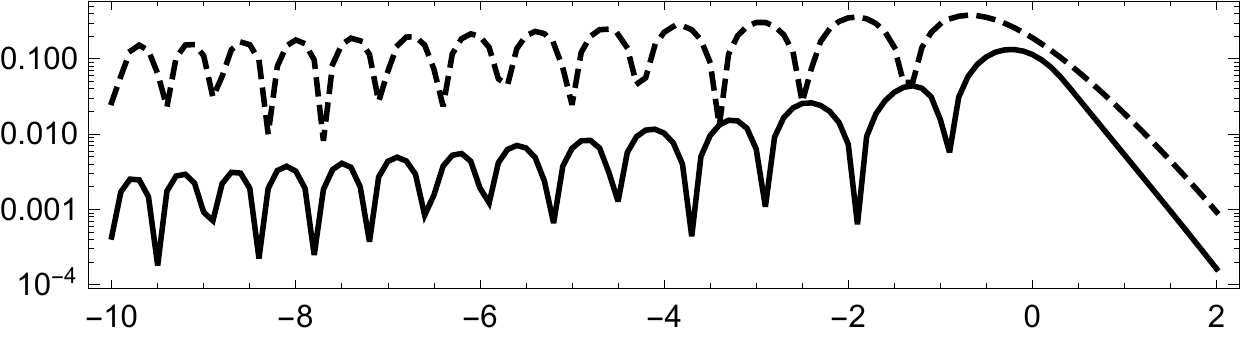}
        \put(197,5){$x$}
    \end{overpic}
    \begin{overpic}[width=0.56\textwidth]{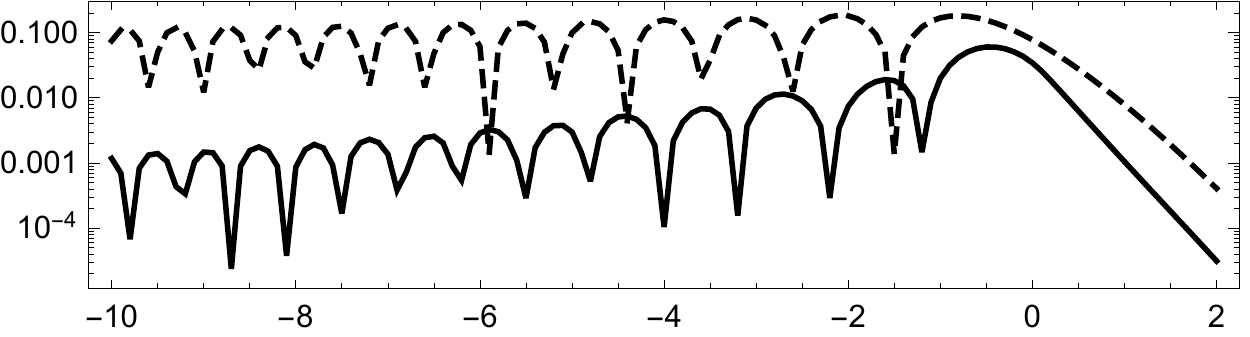}
        \put(197,5){$x$}
    \end{overpic}
    \begin{overpic}[width=0.56\textwidth]{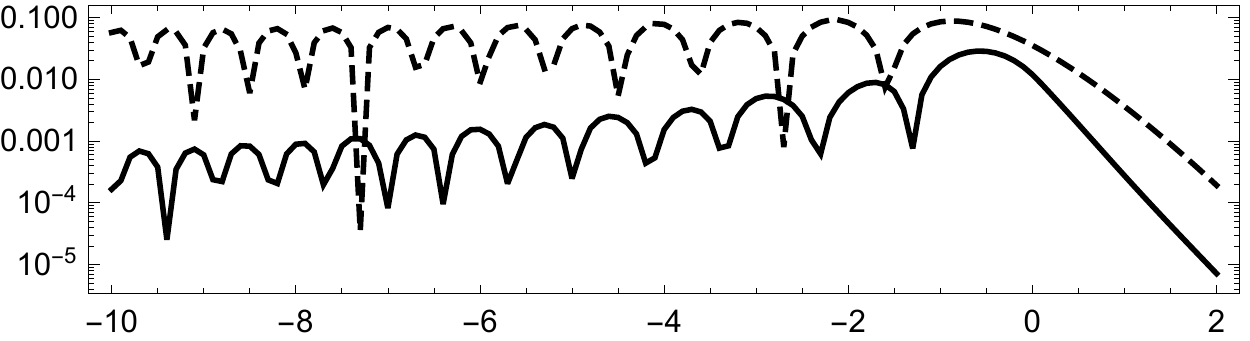}
        \put(197,5){$x$}
    \end{overpic}
\end{center}
\caption{\label{f:compare} A comparison of $E_{\mathrm{model}}(x,0.1)$ (solid) and $E_{\mathrm{LKdV}}(x,.1)$ (dashed) for $a = 1$ (top), $a = 1/2$ (middle) and $a = 1/4$ (bottom).  It is clear that the error $E_{\mathrm{LKdV}}$ dominates the error $E_{\mathrm{model}}$.  The function $u_{\mathrm{KdV}}$ is computed with the method in \cite{Bilman2018}.}
\end{figure}

Our final comparison is for a different take on the small-amplitude limit. Namely, it is well-known that one can approximate the solution of the fully nonlinear KdV equation ($F(u;x,t) = u(x,t)$) with the solution of the linear equation as the amplitude of the initial datum decreases. Similarly, in this limit the solution of our proposed equation \eqref{eq:model} will also approximate the solution of the KdV equation.  The error in such approximation is $O(a^{-1})$ for both cases, but we now provide evidence in Figure~\ref{f:compare} that the constant is much smaller for our model \eqref{eq:model}, and the error decays for large values of $|x|$.  The comparison works as follows.  Consider three solutions of different equations all with initial data given by \eqref{eq:stepdata}:
\begin{enumerate}
    \item $u(x,t) = u(x,t;a,4a)$, the solution of \eqref{eq:model}.
    \item $u_{\mathrm{KdV}}(x,t)$, the solution of the KdV equation, $F(u;x,t) = u(x,t)$.
    \item $u_{\mathrm{LKdV}}(x,t)$, the solution of the linearized KdV equation, $F(u;x,t) = 0$.
\end{enumerate}
We are interested in the two quantities
\begin{align*}
    E_{\mathrm{model}}(x,t) &: = \frac{1}{a}\left| u_{\mathrm{KdV}}(x,t) - u(x,t) \right|,\\
    E_{\mathrm{LKdV}}(x,t) &: = \frac{1}{a}\left| u_{\mathrm{KdV}}(x,t) - u_{\mathrm{LKdV}}(x,t) \right|,
\end{align*}
in the limit $a \downarrow 0$.

\section{Conclusion}

In the present work we have shown that a particular linear model captures a number of qualitative features of the nonlinear dispersive shock phenomenon. The key insight here is to linearize the equation about the initial condition. This leads to a linear differential equation with piecewise constant coefficients. The resulting boundary-value problem is amenable to unified transform techniques.

Our main result is the construction of the linear model and the comparison with the original KdV equation. We emphasize that the unified transform method is only a tool, albeit a very useful one, to solve the proposed linear model. Our results are not contingent on the unified transform method itself. The linear model \eqref{eq:model} reproduces more qualitative features of the nonlinear problem than a naive linearization. Moreover, it is a better short-time linearization. One could in principle, linearize about the solution at some positive time $t$ and repeat the interface problem. Such a procedure may lead to numerical schemes that handle dispersion better. To do so however, would require a `fast' solution method for multiple interfaces. It is possible that the solution presented here is fully expressible in terms of known functions, a combination of re-scaled Airy functions and a stationary solution or perhaps an altogether new special function. These issues will be addressed in future works.

The idea of linearizing about the initial condition also has broader theoretical interest. By moving the initial condition into the differential equation itself, and obtaining a piecewise-constant differential operator, we have realized a simple version of microlocal analysis. In combination with UTM, the linearization procedure presented here can lead to useful (and cheap!) estimates on short-term dynamics. These estimates could serve in establishing local existence of nonlinear dispersive PDEs.

We end on a cautionary note. The faithfulness of our linear model \eqref{eq:model} to reproduce qualitative features of dispersive shocks serves as a warning when comparing theoretical predictions of nonlinear PDEs to experimental results. A fair comparison requires sufficiently long times to allow the nonlinearity to manifest itself in the dynamics.

\section*{Acknowledgements} This research was supported in part by the International Centre for Theoretical Sciences (ICTS) during a visit for participating in the program - Integrable systems in Mathematics, Condensed Matter and Statistical Physics (Code: ICTS/integrability2018/07) and NSF DMS-1945652, DMS-1753185 (TT).

\bibliographystyle{amsplain}
\bibliography{references.bib}

\end{document}